\documentclass{amsart}

\usepackage{amsmath,amssymb,amsthm}

\newtheorem{theorem}{Theorem}[section]
\newtheorem{corollary}[theorem]{Corollary}
\newtheorem{proposition}[theorem]{Proposition}

\theoremstyle{remark}
\newtheorem{remark}{Remark}[section]

\theoremstyle{definition}
\newtheorem{definition}{Definition}[section]

\begin{document}

\title[Absolute Extrema of Invariant Optimal Control Problems]{%
Absolute Extrema of Invariant\\ Optimal Control
Problems\footnote{{\footnotesize Research Report CM06/I-28,
University of Aveiro, 2006. Accepted for publication in the
journal \emph{Commun. Appl. Anal.} (13-Aug-2006).}}}

\author{Cristiana J. Silva}
\author{Delfim F. M. Torres}
\address{Department of Mathematics, University of Aveiro, 3810-193 Aveiro, Portugal}
\email{\texttt{\{cjoaosilva,delfim\}@mat.ua.pt}}
\urladdr{http://www.mat.ua.pt/delfim}

\subjclass[2000]{49J15, 49M30}

\keywords{Optimal control, direct method, symmetries, global extrema}

\maketitle


\begin{abstract}
Optimal control problems are usually addressed with the help of
the famous Pontryagin Maximum Principle (PMP) which gives a
generalization of the classical Euler-Lagrange and Weierstrass
necessary optimality conditions of the calculus of variations.
Success in applying the PMP permits to obtain candidates for a
local minimum. In 1967 a direct method, which permits to obtain
global minimizers directly, without using necessary conditions,
was introduced by Leitmann. Leitmann's approach is connected, as
showed by Carlson in 2002, with ``Carath{\'e}odory's royal road of
the Calculus of variations''. Here we propose a related but
different direct approach to problems of the calculus of
variations and optimal control, which permit to obtain global
minima directly, without recourse to needle variations and
necessary conditions. Our method is inspired by the classical
Noether's theorem and its recent extensions to optimal control. We
make use of the variational symmetries of the problem, considering
parameter-invariance transformations and substituting the original
problem by a parameter-family of optimal control problems.
Parameters are then fixed in order to make the problem trivial, in
some sense. Finally, by applying the inverse of the chosen
invariance-transformation, we get the global minimizer for the
original problem. The proposed method is illustrated, by solving
concrete problems, and compared with Leitmann's approach.
\end{abstract}


\section{Introduction}

The main goal in optimal control is to find a global (or local)
minimizer. One of the most important tools is given by the famous
Pontryagin Maximum Principle (PMP) which is a first order
necessary optimality condition \cite{Pontryagin1962}. The PMP
provides a generalization of the classical Euler-Lagrange and
Weierstrass necessary conditions of the calculus of variations and
permits to obtain candidates for a local minimum. Further analysis
is then needed, to effectively find the extremum.

In 1967 a different approach, based on a coordinate
transformation, was introduced by Leitmann \cite{Leitmann3},
allowing the direct global extremization of certain functional
integrals of the calculus of variations, without the use of
variational methods or field techniques \cite{Leitmann1}. The
method is also valid for multiple integrals of the calculus of
variations \cite{CarlsonLeitmann} and is proved \cite{Carlson} to
be connected  with ``Carath{\'e}odory's royal road of the calculus
of variations'' \cite{Caratheodory}. Here we provide a new look to
Leitmann's approach.

We propose a different direct approach to certain problems of the
calculus of variations and optimal control, which permit to obtain
global minima directly, without recourse to needle variations and
necessary conditions. Differently from Leitmann, our method is
based on the variational symmetries of the problem: a notion
introduced by Emmy Noether in the classical context of the
calculus variations \cite{Noether1} and then extended to the more
general context of optimal control \cite{ejc,Torres4}. Our method
proceeds in three steps: (i) we consider parameter-invariance
transformations of the problem, generalizing the original problem
to an equivalent one; (ii) parameters are then fixed in order to
make the generalized problem trivial in some sense; (iii) finally,
the desired global minimizer is obtained by applying the inverse
of the chosen invariance-transformation and imposing the
fulfilment of the boundary conditions.

The paper is organized as follows. In \S\ref{secPrel} we formulate
the optimal control problem, providing all the necessary
background. In \S\ref{secLeitmethod} we recall Leitmann's approach
and apply it to a simple problem of the calculus of variations.
The same problem is then solved in \S\ref{secDirectMeth}, for
comparison and motivational purposes, by our direct optimization
method. After summarizing the main ideas and steps of the proposed
method, we end \S\ref{secDirectMeth} by considering the minimum
fuel rendezvous of a constant-power rocket. Finally, some
conclusions are presented in \S\ref{sec:Conclusions}.


\section{Preliminaries}
\label{secPrel}

Without loss of generality, we consider the problem of optimal
control in Lagrange form: to minimize an integral functional
\begin{equation}
\label{probCO1} I\left[x(\cdot), u(\cdot)\right] = \int_a^b L
\left(t, x(t), u(t) \right) dt
\end{equation}
subject to a control system
\begin{equation}
\label{controlsystCO1} \dot{x}(t) = \varphi \left(t, x(t),
u(t)\right) \quad \text{a.e. on } [a, b] \, ,
\end{equation}
together with appropriate boundary conditions $x(a) = \alpha$,
$x(b) = \beta$. The Lagrangian $L(\cdot, \cdot, \cdot)$ is a real
function, assumed to be continuously differentiable in $[a,b]
\times \mathbb{R}^{n} \times \mathbb{R}^{m}$; $t \in \mathbb{R}$
is the independent variable; $x:[a,b] \rightarrow \mathbb{R}^{n}$
the vector of state variables; $u:[a,b] \rightarrow \Omega \subset
\mathbb{R}^{m}$ the vector of controls, assumed to be a piecewise
continuous function; and $\varphi :[a,b] \times \mathbb{R}^{n}
\times \mathbb{R}^{m} \rightarrow \mathbb{R}^{n}$ the velocity
vector, assumed to be a continuously differentiable vector
function. In the particular case $\varphi(t,x,u) = u$ one gets the
fundamental problem of the calculus of variations.

The essential concept we are going to use is that of
\emph{equivalence} between two problems of optimal control. In
Carath{\'e}odory's terminology two problems of the calculus of
variations are said to be equivalent when the respective
Lagrangians differ by a total derivative \cite{Caratheodory}. The
importance of this equivalence concept owes to the fact that it
implies the Euler-Lagrange equations to be identical for both
problems. In \cite{Torres6} the following consequence is explored:
\emph{two Carath{\'e}odory-equivalent problems have the same
conservation laws}. It turns out, has proved by E. Noether
\cite{Noether1,Torres4}, that conservation laws are a consequence
of the existence of invariance-transformations (variational
symmetries). The method we propose here is based on the following
trivial remark: the invariance-transformations define a direct
relation between admissible state-control pairs, being
straightforward, from the transformations which define the
equivalence, to obtain a solution for any of the equivalent
problems known the solution for one of them. The variational
symmetries may be found with the help of a computer algebra system
\cite{comPauloLituania05} and, roughly speaking, a given problem
\eqref{probCO1}--\eqref{controlsystCO1} is solved if it admits an
enough rich set of variational symmetries and there exists an
equivalent formulation of the problem with a trivial solution.
This will be illustrated in \S\ref{secDirectMeth}. Now we recall
the notion of invariance (variational symmetry) of an optimal
control problem with respect to a $s$-parameter family of
transformations.

\begin{definition}[\textrm{cf.} \cite{ejc,Torres4}]
\label{definv} Let $h^s(\cdot,\cdot,\cdot)$ be a one-parameter
family of $C^1$ mappings satisfying:
\begin{equation*}
\begin{split}
&h^s:[a,b]\times \mathbb{R}^{n} \times \Omega \longrightarrow
\mathbb{R} \times \mathbb{R}^{n} \times \mathbb{R}^{m} \, , \\
& h^s(t,x,u) = \left( t^s(t,x,u), x^s(t,x,u),
u^s(t,x,u) \right) \, , \\
& h^0(t,x,u)= (t,x,u) \, , \quad \forall (t,x,u) \in [a,b] \times
\mathbb{R}^{n} \times \Omega \, .
\end{split}
\end{equation*}
If there exists a function $\Phi^s(t,x,u) \in C^1\left([a,b],
\mathbb{R}^{n}, \Omega; \mathbb{R} \right)$ such that
\begin{equation}
\label{eq:inv:L} L \circ h^s(t, x(t), u(t)) \frac{d}{dt}t^s
\left(t, x(t), u(t)\right) = L\left(t, x(t), u(t)\right) +
\frac{d}{dt}\Phi^s\left(t, x(t), u(t)\right)
\end{equation}
and
\begin{equation}
\label{eq:inv:CS} \frac{d}{dt}x^s\left(t, x(t), u(t)\right) =
\varphi \circ h^s\left(t, x(t), u(t) \right) \frac{d}{dt}
t^s\left(t, x(t), u(t)\right)
\end{equation}
for all admissible pairs $\left(x(\cdot), u(\cdot)\right)$, then
the optimal control problem
\eqref{probCO1}--\eqref{controlsystCO1} is said to be
\emph{invariant} under the transformations $h^s(t,x,u)$ up to
$\Phi^s(t,x,u)$.
\end{definition}

A parameter-transformation $h^s(t,x,u)$ satisfying all the
conditions of Definition~\ref{definv} is said to be a
\emph{variational symmetry} of the optimal control problem
\eqref{probCO1}--\eqref{controlsystCO1}. These
invariance-transformations are the starting point to our direct
optimization method \S\ref{secDirectMeth}. Next, we review, comment
and illustrate Leitmann's approach.


\section{Remarks on Leitmann's direct optimization method}
\label{secLeitmethod}

G. Leitmann has proposed in 1967 a direct optimization method for
a certain class of scalar problems of the calculus of variations
\cite{Leitmann3}. Leitmann's method is based on the use of
transformations that satisfy a certain functional identity and
permit to obtain, in some cases, absolute extremals directly,
without using variational methods. Since the pioneering work
\cite{Leitmann3} Leitmann has worked on several generalizations of
his method, extending the class of problems to which the method
may be applied: to problems of the calculus of variations where
the trajectory is vector-valued, \textrm{i.e.} $x(t) \in
\mathbb{R}^{n}$, and to problems with side differential conditions
that arise in optimal control \cite{Leitmann1};
by allowing constraints in the form of differential equations and
by considering infinite-horizon problems \cite{Leitmann2}. More
recently, Carlson and Leitmann extended the method to free
problems of the calculus of variations with multiple integrals
\cite{CarlsonLeitmann}. In this section we synthesize Leitmann's
method
\cite{Carlson,CarlsonLeitmann,Leitmann3,Leitmann1,Leitmann2}.
Then, we apply it to solve a simple problem of the calculus of
variations which is used in \S\ref{secDirectMeth} to motivate our
method.


\subsection{Leitmann's main results}

Consider the fundamental problem of the calculus of variations:
\begin{equation}
\label{probasicoCV} I[x(\cdot)]= \int_a^b L\left( t, x(t), \dot{x}
(t)\right) dt \longrightarrow \min \, ,
\end{equation}
where $\dot{x}(t) = \frac{dx(t)}{dt}$, $[a, b]$ is a given fixed
interval, the Lagrangian $L(\cdot, \cdot, \cdot)$ is a real
continuously differentiable function in $[a,b] \times \mathbb{R}
\times \mathbb{R}$, the admissible functions $x(\cdot)$ belong to
$PC^1$ and must satisfy the boundary conditions
\begin{equation}
\label{condfrontCV} x(a)= \alpha  \, ,  \quad x(b) = \beta \, .
\end{equation}

\begin{theorem}
\label{lema1} Let $x=z(t, \tilde{x})$ be a transformation having
an unique inverse $\tilde{x} = \tilde{z}(t, x)$ for $t \in [a,b]$,
such that there is a one-to-one correspondence
\begin{equation*}
x(t) \Leftrightarrow \tilde{x} (t) \, ,
\end{equation*}
for all functions $x(\cdot): [a,b] \rightarrow \mathbb{R}$ in the
class $PC^1$ satisfying \eqref{condfrontCV} and all functions
$\tilde{x}(\cdot): [a,b] \rightarrow \mathbb{R}$ in the class
$PC^1$ satisfying
\begin{equation}
\label{condfrontLema} \tilde{x}(a)=\tilde{z}(a, \alpha)  \, ,
\quad \tilde{x}(b) =\tilde{z}( b, \beta) \, .
\end{equation}
If the transformation $x=z(t, \tilde{x})$ is such that there
exists a function $G:[a,b] \times \mathbb{R} \rightarrow
\mathbb{R}$ such that the functional identity
\begin{equation}
\label{identfuncional} L\left(t, x(t), \dot{x}(t) \right) -
L\left(t, \tilde{x}(t), \dot{\tilde{x}}(t) \right) =
\frac{dG}{dt}{\left(t, \tilde{x}(t) \right)}
\end{equation}
holds, then if $\tilde{x}^{*}(\cdot)$ yields the extremum of
$I[\cdot]$ with $\tilde{x}^{*}(\cdot)$ satisfying
\eqref{condfrontLema}, $x^*(t)=z\left( t, \tilde{x}^*(t)\right)$
yields the extremum of $I[\cdot]$ for $x^*(\cdot)$ satisfying
\eqref{condfrontCV}.
\end{theorem}

\begin{remark}
\label{eq:obs} To the best of our knowledge, no one has
interpreted \eqref{identfuncional} before as being Noether's
invariance condition \eqref{eq:inv:L} in the particular case where
no transformation of time is considered, \textrm{i.e.} $t^s = t$.
Instead of \eqref{identfuncional}, the method we propose here is
based on the more rich set of identities
\eqref{eq:inv:L}-\eqref{eq:inv:CS}.
\end{remark}

There is a one-to-one correspondence between the minimizers of
problem \eqref{probasicoCV}-\eqref{condfrontCV} and the minimizers
of the integral functional $I[\tilde{x}(\cdot)]= \int_a^b
L\left(t, \tilde{x}(t), \dot{\tilde{x}}(t) \right) dt$ in the
class of functions $\tilde{x}(\cdot) \in PC^1$ satisfying the
boundary conditions \eqref{condfrontLema}. Moreover, the
transformation $x=z(t,\tilde{x})$ and its inverse $\tilde{x} =
\tilde{z}(t, x)$ give us the desired correspondence.

\begin{corollary}
\label{corol1} For the validity of Theorem \ref{lema1}, the
Lagrangian $L(\cdot, \cdot, \cdot )$, together with the
transformation $x=z(t, \tilde{x})$, must be such that the
left-hand side of the functional identity \eqref{identfuncional}
is linear with respect to $\dot{\tilde{x}}(t)$.
\end{corollary}

The main difficulty in applying Leitmann's method
(Theorem~\ref{lema1}) resides in finding the admissible
transformations $x=z(t, \tilde{x})$. Leitmann has restricted
himself to two situations for which it is easy to find the
admissible transformations: (i) Corollary~\ref{corol1} is
trivially satisfied if $L(\cdot, \cdot, \cdot)$ is linear in its
third argument; (ii) it can also be readily satisfied for
$L(\cdot, \cdot, \cdot)$ quadratic in its third argument,
\textrm{i.e.} for $L(\cdot, \cdot, \cdot)$ of the form
\begin{equation}
\label{formaquad} L(t, x, p) = a(t)p^2 + b(t,x)p + c(t,x) \, ,
\end{equation}
with $a(t) \neq 0$ for $t \in [a,b]$.

\begin{corollary}
\label{corol2} For a Lagrangian of type \eqref{formaquad} the
class of admissible transformations that satisfy
Corollary~\ref{corol1} is of the form $x = z(t, \tilde{x}) = \pm
\tilde{x}+ f(t)$.
\end{corollary}

\begin{remark}
Using our Remark~\ref{eq:obs} we can benefit of a well-developed theory
\cite{PauloBrasil,comPauloLituania05,GouveiaTorresRochaPolonia05}
on how to find Noether's invariance transformations, without the
need to restrict ourselves to Lagrangians which are linear or
quadratic in the velocity.\footnote{A computer algebra package
to compute variational symmetries,
by Paulo D. F. Gouveia and Delfim F. M. Torres,
is available from the \emph{Maple Application Center}:
\texttt{http://www.maplesoft.com/applications/app\_center\_view.aspx?AID=1983}}
Therefore, the method we propose is
applicable to a more wide class of optimization problems.
\end{remark}


\subsection{An example}

Let us apply Leitmann's method to the following simple problem of
optimal control ($a < b$):
\begin{equation}
\label{probasCO}
\begin{gathered}
I[u(\cdot)]= \int_a^b \left(u(t) \right)^2 dt \longrightarrow \min \, ,\\
\dot{x}(t) = u(t) \, , \\
x(a)= \alpha \, , \quad x(b) = \beta \, .
\end{gathered}
\end{equation}
This is a problem \eqref{probCO1}-\eqref{controlsystCO1} with
$\varphi = u$, so we can write \eqref{probasCO} as a problem
\eqref{probasicoCV}-\eqref{condfrontCV} of the calculus of
variations:
\begin{equation*}
\begin{gathered}I[x(\cdot)]= \int_a^b \left(\dot{x}(t)\right)^2 dt \longrightarrow \min
\, , \quad x(a)= \alpha \, , \quad x(b) = \beta \, .
\end{gathered}
\end{equation*}
The Lagrangian $L(t,x,p) = p^2$ is of type \eqref{formaquad},
thus, by Corollary \ref{corol2}, the class of admissible
transformations of Theorem~\ref{lema1} has the form $x= z(t,
\tilde{x})= \pm \tilde{x} + f(t)$, where $f(t)$ is some
differentiable function. We consider, without loss of generality,
the transformation $x= z(t, \tilde{x})= \tilde{x} + f(t)$. Then,
\begin{equation*}
L\left(t, f(t)+ \tilde{x}, f'(t) + \tilde{p} \right) - L\left(t,
\tilde{x}, \tilde{p}' \right) = \left( f'(t) \right)^2 +
2f'(t)\tilde{p} \, ,
\end{equation*}
and from the functional identity \eqref{identfuncional} we get
\begin{equation*}
\frac{\partial G}{\partial t}(t, \tilde{x}) = \left( f'(t)
\right)^2 \, , \quad \frac{\partial G}{\partial \tilde{x}}(t,
\tilde{x}) = 2 f'(t) \, .
\end{equation*}
On the other hand,
\begin{equation*}
\frac{\partial^2 G}{ \partial \tilde{x} \partial t}(t, \tilde{x})
= \frac{\partial^2 G}{\partial t \partial \tilde{x}}(t, \tilde{x})
\, ,
\end{equation*}
and we conclude that
\begin{equation}
\label{eqELag}
 2f''(t) =0  \, ,
\end{equation}
that is,
\begin{equation}
\label{ft1} f(t)=c_1 + c_2 t \, ,
\end{equation}
with $c_1$ and $c_2$ constants. We now determine function $G(t,
\tilde{x})$. Substituting \eqref{ft1} into $\frac{\partial
G}{\partial \tilde{x}}(t, \tilde{x}) = 2 f'(t)$ it follows that
$\frac{\partial G}{\partial \tilde{x}}(t, \tilde{x}) = 2 c_2$, and
integrating with respect to $\tilde{x}$ we arrive to
\begin{equation*}
G(t, \tilde{x}) = \int (2 c_2)d\tilde{x} = 2c_2 \tilde{x} + h(t)
\, ,
\end{equation*}
where $h(t)$ is still to be determined. For that, we differentiate
the last expression with respect to $t$ and compare the result
with $\frac{\partial G}{\partial t}(t, \tilde{x})$:
\begin{equation*}
\frac{\partial}{\partial t}\left( 2c_2 \tilde{x} + h(t) \right) = h'(t) \, ,
\end{equation*}
and since $\frac{\partial G}{\partial t}(t, \tilde{x}) = (c_2)^2$
we must have
\begin{equation*}
h'(t) = (c_2)^2 \Leftrightarrow h(t) = (c_2)^2t + c_3 \, ,
\end{equation*}
where $c_3$ is an arbitrary constant. Therefore,
$G(t, \tilde{x}) = 2c_2 \tilde{x} + (c_2)^2t + c_3$.
We have all the necessary ingredients to apply Theorem~\ref{lema1}. We
consider the trivial problem
\begin{equation*}
I[\tilde{x}(\cdot)]= \int_a^b \left( \dot{\tilde{x}}(t) \right)^2
dt \longrightarrow \min \, , \qquad \tilde{x}(a)=0  \, , \quad
\tilde{x}(b)=0 \, ,
\end{equation*}
which admits the global minimizer $\tilde{x}^*(t) \equiv 0$ (the
original problem \eqref{probasCO} is trivial when $\alpha =
\beta$; we are interested to solve \eqref{probasCO} in the case
$\alpha \ne \beta$). To obtain the solution of problem
\eqref{probasCO} we just need to choose $c_1$ and $c_2$ in
\eqref{ft1} in such a way $f(a) = \alpha$ and $f(b) = \beta$, \textrm{i.e.}
\begin{equation*}
\begin{cases}
f(a)=\alpha \\[0.1cm]
f(b)= \beta
\end{cases}
\Leftrightarrow
\begin{cases}
c_1 + c_2 a=\alpha \\[0.1cm]
c_1 + c_2 b= \beta
\end{cases}
\Leftrightarrow
\begin{cases}
c_1 =\frac{\beta a - b \alpha}{a-b}  \\[0.1cm]
c_2= \frac{\alpha - \beta}{a-b} \, .
\end{cases}
\end{equation*}
The global minimizer for problem \eqref{probasCO} is given by
$x^*(t) = \tilde{x}^*(t) + f(t) = 0 + c_1 + c_2 t$:
\begin{equation}
\label{solLeitm} x^*(t) = \frac{\beta a - b \alpha}{a-b} +
\frac{\alpha - \beta}{a-b} t \, .
\end{equation}

We remark that \eqref{solLeitm} satisfies \eqref{eqELag}, and that
\eqref{eqELag} is nothing more than the Euler-Lagrange equation of
\eqref{probasCO}. However, the Euler-Lagrange equation only gives
a candidate for local minimizer, \textrm{i.e.} we are not sure if
the candidate is indeed a local minimizer. Leitmann's method has
given much more: \eqref{solLeitm} is the global minimizer of
\eqref{probasCO}. Next section gives an alternative direct
optimization method, which we claim to be more broad in
application.


\section{A new direct optimization method}
\label{secDirectMeth}

Our direct optimization method is of simple comprehension and is
applicable to a wider class of optimal control problems. We first
show how it can be applied to problem \eqref{probasCO}.


\subsection{Motivational example}

The initial step of our method is the determination of the
parametric transformations under which the problem is invariant
(see Definition~\ref{definv}). In respect to this, the techniques
found in \cite{comPauloLituania05,Torres4} are useful.

\begin{proposition}\cite[Ex.~1]{ejc}
\label{propInv1} Problem \eqref{probasCO} is invariant
up to $\Phi^s\left(t,x\right) = s^2 t + 2 s x$,
in the sense of Definition~\ref{definv}, under the
$s$-parameter transformations ($s \in \mathbb{R}$)
\begin{equation}
\label{transf}
t^s = t \, , \quad x^{s}= x + s t \, , \quad u^{s} = u + s \, .
\end{equation}
\end{proposition}

\begin{proof}
We begin by showing \eqref{eq:inv:L}:
\begin{equation}
\label{eq:invL:ex}
\begin{split}
\tilde{I} &= \int_a^b \left(u^{s}(t)\right)^2 dt
= \int_a^b \left(u(t)+s \right)^2 dt
= \int_a^b \left(u^2(t)+s^2 + 2su(t) \right) dt \\
&= \int_a^b u^2(t) dt + \int_a^b \left( s^2 + 2su(t) \right) dt
= I + \int_a^b \frac{d}{dt}\left(s^2t + 2sx(t) \right)dt \\
&= I + \Phi^s(b,\beta) - \Phi^s(a,\alpha) \, .
\end{split}
\end{equation}
We remark that the minimizer of $\tilde{I}[\cdot]$ coincide with
the one of $I[\cdot]$: $\Phi^s(a,\alpha)$ and $\Phi^s(b,\beta)$
are constants and adding a constant in the functional does not
change the minimizer. It remains to prove the control invariance
condition \eqref{eq:inv:CS}:
\begin{equation}
\label{eq:invPhi:ex}
\frac{d}{dt}\left(x^{s}(t) \right)= \frac{d}{dt} \left( x(t) + st
\right) = \dot{x}(t) + s = u(t)+ s = u^{s}(t) \, .
\end{equation}
Equalities \eqref{eq:invL:ex} and \eqref{eq:invPhi:ex} prove that
problem \eqref{probasCO} is invariant under the one-parameter
transformations \eqref{transf} up to the gauge term $\Phi^s$.
\end{proof}

Using the invariance transformations \eqref{transf} we generalize
problem \eqref{probasCO} to a parameter family of problems which
include the original problem for $s = 0$: we substitute $x(\cdot)$
and $u(\cdot)$ in \eqref{probasCO} respectively by $x^s(\cdot)$
and $u^s(\cdot)$, obtaining
\begin{equation}
\label{intCO-eq}
\begin{gathered}
I^s[u^s(\cdot)]=
\int_a^b \left(u^{s}(t) \right)^2 dt \longrightarrow \min \, ,\\
\dot{x}^{s}(t)= u^{s}(t) \, , \\
x^{s}(a) = \alpha + sa \, , \quad x^{s}(b) = \beta + sb \, , \quad
s\in \mathbb{R} \, .
\end{gathered}
\end{equation}
Problem \eqref{probasCO} is nontrivial for $\alpha \ne \beta$, but
the crucial point is that there exists always a problem in the
parameter family of problems \eqref{intCO-eq}, \textrm{i.e.} there
exists always a specific value of $s$, which only depend on the
concrete values of $\alpha$, $\beta$, $a$ and $b$, admitting the
trivial global minimizer $u^s(t)=0$ $\forall$ $t \in [a,b]$. The
invariance properties asserted by Proposition~\ref{propInv1} give
the general solution to our original problem \eqref{probasCO}
from the trivial solution of this $s$-chosen problem.

\begin{proposition}
Function \eqref{solLeitm} is a global minimizer of problem
\eqref{probasCO}.
\end{proposition}
\begin{proof}
It is clear that $I^s \ge 0$ and that $I^s = 0$ if $u^s(t)
\equiv 0$. From the control system $\dot{x}^{s}(t)= u^{s}(t)$,
$u^s(t) \equiv 0$ implies that $x^s(a) = x^s(b)$:
\begin{equation*}
\alpha + sa = \beta + sb \Leftrightarrow
s = \frac{\beta - \alpha}{a-b} \, .
\end{equation*}
Hence, the global minimizing trajectory of problem \eqref{intCO-eq}
for $s=\frac{\beta - \alpha}{a-b}$ is given by
\begin{equation*}
x^{s}(t)= \alpha + sa \Leftrightarrow
x^{s}(t)= \frac{\beta a - \alpha b}{a-b} \, .
\end{equation*}
We solve \eqref{probasCO} using the inverse functions
of the variational symmetries \eqref{transf}:
\begin{equation*}
\begin{cases}
u(t) = u^{s}(t) - s \\[0.1cm]
x(t) = x^{s}(t) - st
\end{cases}
\Leftrightarrow
\begin{cases}
u(t) = \frac{\alpha - \beta}{a-b} \\[0.1cm]
x(t) = \frac{\beta a - \alpha b}{a-b} - \frac{\beta - \alpha}{a-b}t \, .
\end{cases}
\end{equation*}
We have just found the global minimizer \eqref{solLeitm}
of problem \eqref{probasCO}.
\end{proof}


\subsection{The method}
\label{subsec:TheMethod}

As just illustrated, our direct optimization method permits to
find global extremizers (minimizers or maximizers) of sufficiently
rich invariant optimal control problems. The method consists of
the following four steps:

\begin{itemize}
\item[(1)] Determine parameter invariant transformations
$t^s$, $x^s_i$, and $u^s_j$, $i=1,\ldots,n$, $j=1,\ldots,m$, under
which the problem is invariant (\textrm{cf.}
Definition~\ref{definv}). The results in
\cite{comPauloLituania05,Torres4} are useful.

\item[(2)] Applying the parameter transformations found in the
previous step, write the generalized problem together with the
generalized boundary conditions, \textrm{i.e.} substitute
$x_i(\cdot)$ and $u_j(\cdot)$ respectively by $x_i^s(\cdot)$ and
$u_j^s(\cdot)$, $i=1, \ldots, n$ and $j=1, \ldots, m$.

\item[(3)] Analyze the generalized problem and
determine a specific value for the parameters for which it is easy
to find a global optimal solution.

\item[(4)] Define the inverse of the transformations $t^s$,
$x^s$, and $u^s$, for the particular choice of parameters $s$
fixed on step $(3)$, and obtain a
global solution to the initial problem.

\end{itemize}

We shall now apply our simple method to the minimum fuel
rendezvous problem of a constant-power rocket.


\subsection{An application}

Let us consider the problem of minimizing the amount of fuel
consumed by a rocket operating at constant propulsive power. This
is a classical problem of optimal control, ``solved'' by the
Pontryagin Maximum Principle in most books (see \textrm{e.g.}
\cite{Leitmann4,MackiStrauss}). We assume the following situation:
(i) a positive prescribed transfer time $\tau
= t_1 - t_0$ is given; (ii) at the end (at time $t_1$) the
rocket car is to be at the origin with zero-velocity; (iii)
the rocket is initially on the negative axis (at a given position
$-\alpha$, $\alpha > 0$). Thus, we have:
\begin{equation}
\label{eq:rocket:Lag}
\begin{gathered}
\int_{t_0}^{t_1} u^2(t) dt \longrightarrow \min \, ,
\quad t_1 - t_0 = \tau \, , \\
\begin{cases}
\dot{x}_1(t)= x_2(t) \, , \\[0.1 cm]
\dot{x}_2(t)= u(t) \, ,
\end{cases} \\
x_1(t_0)= - \alpha \, , \quad x_1(t_1)= 0 \, ,
\quad x_2(t_1) = 0 \, ,
\end{gathered}
\end{equation}
where $t$ is the time variable, $x_1$ the position, $x_2$ the
velocity, and $u$ is the acceleration due to the thrust. We are
assuming that the thrust-acceleration is not constrained,
\textrm{i.e.} $|u(t)| < \infty$, and that $\tau > 0$ and $\alpha > 0$ are
given. The thrust-acceleration program that results in the minimum
fuel consumption can also be obtained by Leitmann's method
(\textrm{cf.} \cite[\S~9]{Leitmann1}) but the analysis is
enough-complex: it is not easy to guess functions $f_1(t)$ and
$f_2(t)$ of Corollary~\ref{corol2}, associated respectively with
$x_1(t)$ and $x_2(t)$. Here we show that there exists a simple way
to obtain a global minimizer to problem \eqref{eq:rocket:Lag}.

\begin{proposition}
A global minimizer of problem \eqref{eq:rocket:Lag} is given by
\begin{equation}
\label{eq:sol:RocCar}
\begin{gathered}
x_1(t)= -\frac{\alpha}{\tau^2} t^2
+ \frac{2 \alpha}{\tau} \left( \frac{t_0}{\tau} + 1 \right)t
- \frac{2 \alpha}{\tau} \left( \frac{t_0}{2\tau} + 1
\right) t_0 - \alpha \, , \\
x_2(t) = - \frac{2 \alpha}{\tau^2} t + \frac{2 \alpha}{\tau}
\left( \frac{t_0}{\tau} + 1 \right) \, , \\
u(t) =  - \frac{2 \alpha}{\tau^2} \, .
\end{gathered}
\end{equation}
\end{proposition}

\begin{proof}
We follow the four-step method of \S\ref{subsec:TheMethod}.

\medskip

(1) Problem \eqref{eq:rocket:Lag} is
invariant under the parameter transformations ($s \in \mathbb{R}$)
\begin{equation}
\label{eq:trf:PR}
t^s = t \, , \ \ x_1^{s}(t)= x_1(t) + \frac{s^2}{2}t^2 \, , \ \
x_2^{s}(t)= x_2(t) + s^2t \, , \ \ u^{s}(t)= u(t) + s^2 \, ,
\end{equation}
up to $\Phi^s\left(t, x_2\right) = s^4 t + 2s^2 x_2$:
the functional is invariant,
\begin{equation*}
\begin{split}
\int_{t_0}^{t_1} \left(u^{s}(t)\right)^2 dt
&= \int_{t_0}^{t_1} \left(u(t)+s^2 \right)^2 dt
= \int_{t_0}^{t_1} \left(u^2(t)+s^4 + 2s^2u(t) \right) dt \\
&= \int_{t_0}^{t_1} u^2(t) dt + \int_{t_0}^{t_1} \left( s^4 + 2s^2u(t) \right) dt \\
&= \int_{t_0}^{t_1} u^2(t) dt + \int_{t_0}^{t_1}
\frac{d}{dt}\left(s^4 t + 2s^2x_2(t) \right)dt \, ;
\end{split}
\end{equation*}
as well as the control system,
\begin{equation*}
\begin{cases}
\frac{d}{dt}\left(x_1^{s}(t) \right)= \frac{d}{dt} \left( x_1(t) +
\frac{s^2}{2}t^2
\right) = \dot{x_1}(t) + s^2t = x_2(t)+ s^2t = x_2^{s}(t) \, , \\[0.1 cm]
\frac{d}{dt}\left(x_2^{s}(t) \right)= \frac{d}{dt} \left( x_2(t) +
s^2t \right) = \dot{x_2}(t) + s^2 = u(t)+ s^2 = u^{s}(t) \, .
\end{cases}
\end{equation*}

(2) The generalized problem takes the following form:
\begin{equation}
\label{eq:Appl:GenPrb}
\begin{gathered}
I^s[u^{s}(\cdot)] = \int_{t_0}^{t_1} \left(u^{s}(t)\right)^2 dt \longrightarrow \min \, , \\
\begin{cases}
\dot{x}_1^s(t)= x_2^s(t) \, , \\[0.1 cm]
\dot{x}_2^s(t)= u^s(t) \, ,
\end{cases} \\
x_1^s(t_0)= \frac{s^2 t_0^2}{2} - \alpha \, , \quad
x_1^s(t_1)= \frac{s^2 t_1^2}{2} \, , \quad
x_2^s(t_1)= s^2 t_1 \, .
\end{gathered}
\end{equation}
For $s = 0$ problem \eqref{eq:Appl:GenPrb}
reduces to \eqref{eq:rocket:Lag}.

\medskip

(3) $I^s \ge 0$ $\forall$ $u^s(\cdot)$,
and $I^s = 0$ if $u^s(t) \equiv 0$.
From the control system we have
for $u^s(t) \equiv 0$ that
$x_2^s(t) = c_1$, and
$x_1^s(t) = c_1 t + c_2$,
where $c_1$ and $c_2$ are constants.
From the generalized boundary condition $x_2^s(t_1)= s^2 t_1$,
it follows that $c_1 = s^2 t_1$.
Then, $x_2^s(t) = s^2 t_1$,
$x_1^s(t) = s^2 t_1 t + c_2$. Using the boundary conditions
for $x_1^s(\cdot)$ we arrive to $c_2 = - \frac{s^2 t_1^2}{2}$
and $s^2 = \frac{2 \alpha}{\tau^2}$. Therefore, a global
minimizer to problem \eqref{eq:Appl:GenPrb} with
$s = \pm \frac{\sqrt{2 \alpha}}{\tau}$ is given by
\begin{equation*}
x_1^s(t) = \frac{2 \alpha}{\tau^2} t_1 t - \frac{\alpha t_1^2}{\tau^2} \, , \quad
x_2^s(t) = \frac{2 \alpha}{\tau^2} t_1 \, , \quad u^s(t) = 0 \, , \quad t \in [t_0,t_1] \, .
\end{equation*}

(4) The global solution to problem \eqref{eq:rocket:Lag}
is obtained using the inverses of transformations \eqref{eq:trf:PR}
for $s^2 = \frac{2 \alpha}{\tau^2}$:
\begin{equation}
\label{eq:solEquiv:RocCar}
\begin{gathered}
x_1(t)= x_1^{s}(t) - \frac{s^2}{2}t^2 = \frac{2\alpha}{\tau^2} t_1 t
- \frac{\alpha t_1^2}{\tau^2} - \frac{\alpha}{\tau^2} t^2 \, , \\
x_2(t) = x_2^{s}(t) - s^2t = \frac{2 \alpha}{\tau^2} t_1 - \frac{2 \alpha}{\tau^2} t \, ,\\
u(t) = u^{s}(t) - s^2 = - \frac{2 \alpha}{\tau^2} \, .
\end{gathered}
\end{equation}
It is a simple exercise to see that \eqref{eq:sol:RocCar}
and \eqref{eq:solEquiv:RocCar} are equivalent.
\end{proof}


\section{Conclusions}
\label{sec:Conclusions}

In the calculus of variations, as well as in the more general
setting of optimal control, the problem of minimizing an integral
functional is the main issue, in general a difficult one. The
standard way to attack such problems relies on necessary
optimality conditions, which give candidates for a local minimum.
A direct method for addressing some problems of the calculus of
variations which are linear or quadratic in velocity (control) was
introduced by Leitmann and further improved by Carlson, providing
global minimizers directly, without using necessary conditions.
Here we propose a different, simpler, and more wide applicable
direct method for problems of optimal control: (i) different
because instead of using transformations which keep the problem
invariant in Carath{\'e}odory's sense, as in the method of
Leitmann-Carlson, our method is based on transformations which
keep the problems invariant in Noether's sense; (ii) simpler in
finding the admissible transformations; (iii) more general because
it easily covers Lagrangians which are not linear or quadratic in
the control variables.


\section*{Acknowledgments}

This work was partially supported by the
\emph{Portuguese Foundation for Science and Technology} (FCT),
cofinanced by the European Community Fund FEDER/POCTI,
through the \emph{Control Theory Group} (cotg) of the
\emph{Centre for Research on Optimization and Control} (CEOC):
\texttt{http://ceoc.mat.ua.pt}



\end{document}